\documentclass[11pt]{article}
\usepackage{amsthm,amsmath,latexsym,amssymb}

\newcommand{\bP}{{\rm |\kern-.15em P}}
\newcommand{\Q}{\kern.3em\rule{.07em}{.65em}\kern-.3em{\rm Q}}
\newcommand{\R}{{\rm I\kern-.15em R}}
\newcommand{\D}{{\rm |\kern-.15em D}}
\newcommand{\h}{{\rm |\kern-.15em H}}
\newcommand{\C}{\kern.3em\rule{.07em}{.65em}\kern-.3em{\rm C}}
\newcommand{\T}{{\rm T\kern-.35em T}}

\theoremstyle{plain}
\newtheorem{theorem}{Theorem}[section]

\newtheorem{corollary}[theorem]{Corollary}

\theoremstyle{definition}
\newtheorem{definition}[theorem]{Definition}

\theoremstyle{remark}
\newtheorem{remark}[theorem]{Remark}

\begin{document}
\title{On an integral equation of Lieb}
\author{Ronen Peretz}
 
\maketitle

\begin{abstract}
We prove that the weakly singular, non-linear convolution integral equation $\int_{\mathbb{R}^n}|x-y|^{-\lambda}f(y)dy=f(x)^{p-1}$,
where $0<\lambda<n$, and $p=2n/(2n-\lambda)$ has at least two non-equivalent solutions. This answers a problem of Elliott Lieb.
We also prove certain orthogonality relations among linear differential forms with constant coefficients related to the corresponding 
type of convolution operators. Finally, we discuss the regularity of the solutions of such non-linear integral equations over not necessarily
bounded open subsets of $\mathbb{R}^n$.
\end{abstract}

\section{Introduction of the results}
We divide the results in this paper into three parts. In {\bf the first part} we note that a certain non-linear integral equation introduced by Elliott
Lieb has at least two essentially different solutions. One of the solutions has an isolated singular point where the function tends to infinity.

{\bf The second part} presents a multitude of integral identities each connecting two solutions of
Lieb's equation. These integral identities contain the images of the two solutions under linear differential forms with constant coefficients. They
point out to two properties connecting any two such solutions. One property is an orthogonality property and the second property is a certain commutativity. 

Finally, in {\bf the third part} we present regularity results of the solutions
to the Lieb integral equation. The bottom line is that except for isolated singularities were the solutions tend to infinity, they are, in fact
smooth functions elsewhere. One can compute effectively their smoothness degree. 

{\bf (I)} In \cite{Lieb} Elliott H. Lieb computes the sharp constants for certain parameter values in the Hardy-Littlewood-Sobolev inequality. The paper
also deals with related inequalities: the Sobolev inequality, doubly weighted Hardy-Littlewood-Sobolev inequality and the weighted Young inequality
(as A. Sokal called it). Theorem 3.1 on page 359 computes the unique maximizing function and the sharp constant for certain values of the
parameters of the first inequality. The maximizing functions should satisfy the following integral equation
\begin{equation}
\label{eq1}
\int_{\mathbb{R}^n}\left|x-y\right|^{-\lambda}f(y)dy=f(x)^{p-1},\,\,\,0<\lambda<n,\,\,\,p=\frac{2n}{2n-\lambda}.
\end{equation}
However, as claimed on page 361 of \cite{Lieb} "We do not that (3.9)" (the above (\ref{eq1}) equation) "has an (essentially) unique solution-even if we restrict
to the SSD category-and we shall offer no proof of this kind of uniqueness. {\it This is an open problem!}". Indeed there is no uniqueness for we have the
following:

\begin{theorem}
The following function is a solution of equation (\ref{eq1}):
\begin{equation}
\label{eq2}
f(x)=C(n,\lambda)|x|^{-(n-\lambda/2)},
\end{equation}
where:
\begin{equation}
\label{eq3}
C(n,\lambda)=
\end{equation}
$$
=\left\{\pi^{n/2}\left(\Gamma\left(\frac{n}{2}-\frac{\lambda}{2}\right)\Gamma\left(\frac{\lambda}{4}\right)^2\right)\left/
\left(\Gamma\left(\frac{\lambda}{2}\right)\Gamma\left(\frac{n}{2}-\frac{\lambda}{4}\right)^2\right)\right.\right\}^{-(2n-\lambda)/(2(n-\lambda))}.
$$
\end{theorem}

\begin{remark}
Lieb proved that his equation (\ref{eq1}) has the following solution (which is unique among the maximizing functions of the inequality he was considering):
\begin{equation}
\label{eq4}
f_L(x)=L(n,\lambda)(1+|x|^2)^{-n/p}=L(n,\lambda)(1+|x|^2)^{-(n-\lambda/2)}.
\end{equation}
Here $L(n,\lambda)$ is a constant depending on the dimension and on the parameter $\lambda$. We note that the solution in Theorem 1.1 is singular
at the origin and is not a reflection of Lieb's solution. It is certainly not a conformal image of it ($f_L$ is non-singular and bounded). Thus
Lieb's equation (\ref{eq1}) exhibits at least two non-equivalent solutions.
\end{remark}

Based on the coming results we can deduce many integral identities. Here is an example:

\begin{corollary}
$$
L(n,\lambda)^{\lambda/(2n-\lambda)}C(n,\lambda)\int_{\mathbb{R}^n}|x|^{-(n-\lambda/2)}(1+|x|^2)^{-\lambda/2}dx=
$$
\begin{equation}
\label{eq5}
\end{equation}
$$
=L(n,\lambda)C(n,\lambda)^{\lambda/(2n-\lambda)}\int_{\mathbb{R}^n}|x|^{-\lambda/2}(1+|x|^2)^{-(n-\lambda/2)}dx.
$$
\end{corollary}
We note the symmetric relations between the two solutions within the last identity. It is a kind of commutativity that takes the power $p-1$ into
an account.

{\bf (II)} The identity in the last corollary is the zero'th integral identity out of infinitely many possible integral identities that connect two solutions of the Lieb integral equation (\ref{eq1}).
\begin{definition}
For $\alpha=(\alpha_1,\ldots,\alpha_n)\in (\mathbb{Z}^+\cup\{0\})^n$ we denote:
$$
D_{\alpha}^{(x)}(h(x))=\frac{\partial^{|\alpha|}h(x)}{\partial x_1^{\alpha_1}\ldots\partial x_n^{\alpha_n}},\,\,\Lambda^{(x)}=\sum_{\alpha} a_{\alpha}
D_{\alpha}^{(x)}\,\,\,{\rm where}\,\,a_{\alpha}\in\mathbb{R}.
$$
\end{definition}

Here are some integral relations between solutions of equation (\ref{eq1}). Obvious generalizations hold true for other similar kernels. These include some
orthogonality relations and some commutativity relations:

\begin{theorem}
Let $f(x)$ and $g(x)$ be solutions of equation (\ref{eq1}). Then assuming the convergence of the integrals below, $\forall\,\alpha,\beta\in(\mathbb{Z}^+\cup\{0\})^n$ we have:
\begin{equation}
\label{eq6}
\int_{\mathbb{R}^n}D_{\beta}^{(x)}(g(x))D_{\alpha}^{(x)}(f(x)^{p-1})dx=\int_{\mathbb{R}^n}D_{\alpha}^{(x)}(f(x))D_{\beta}^{(x)}(g(x)^{p-1})dx,
\end{equation}
$$
(-1)^{|\beta|}\int_{\mathbb{R}^n}D_{\beta}^{(x)}(f(x))D_{\alpha}^{(x)}(f(x)^{p-1})dx=
$$
\begin{equation}
\label{eq7}
\end{equation}
$$
=(-1)^{|\alpha|}\int_{\mathbb{R}^n}D_{\alpha}^{(x)}(f(x))D_{\beta}^{(x)}(f(x)^{p-1})dx,
$$
$$
(-1)^{|\alpha|+|\beta|}=-1
$$
\begin{equation}
\label{eq8}
\Downarrow
\end{equation}
$$
\int_{\mathbb{R}^n}D_{\beta}^{(x)}(f(x))D_{\alpha}^{(x)}(f(x)^{p-1})dx=\int_{\mathbb{R}^n}D_{\alpha}^{(x)}(f(x))D_{\beta}^{(x)}(f(x)^{p-1})dx=0,
$$
\end{theorem}

So it is natural to make the following:
\begin{definition}
$$
E=\{\sum_{\alpha\in I}a_{\alpha}D_{\alpha}^{(x)}(\cdot)\,|\,\forall\,\alpha\in I,\,(-1)^{|\alpha|}=1,\,a_{\alpha}\in\mathbb{R}\},
$$
$$
O=\{\sum_{\beta\in J}b_{\beta}D_{\beta}^{(x)}(\cdot)\,|\,\forall\,\beta\in J,\,(-1)^{|\beta|}=-1,\,b_{\beta}\in\mathbb{R}\}.
$$
\end{definition}

\begin{theorem}
Let $f(x)$ and $g(x)$ be solutions of equation (\ref{eq1}). Suppose that $\Lambda=\Lambda_e+\Lambda_o$, $\Omega=\Omega_e+\Omega_o$, 
where $\Lambda_e,\Omega_e\in E$, and $\Lambda_o,\Omega_o\in O$. Then: \\
(a) Assuming convergence of the integrals below we have:
\begin{equation}
\label{eq9}
\end{equation}
$$
\int_{\mathbb{R}^n}\Lambda(f)\Omega(g^{p-1})dx=\int_{\mathbb{R}^n}\Lambda(f^{p-1})\Omega(g)dx,
$$
$$
\int_{\mathbb{R}^n}\Lambda_{e}(f(x))\Lambda_{o}(f(x)^{p-1})dx=\int_{\mathbb{R}^n}\Lambda_{e}(f(x)^{p-1})\Lambda_{o}(f(x))dx=0.
$$
(b) Assuming the convergence of the integrals below we have:
\begin{equation}
\label{eq10}
\end{equation}
$$
\int_{\mathbb{R}^n}\Lambda(f)\Omega(f^{p-1})dx=\int_{\mathbb{R}^n}\Lambda(f^{p-1})\Omega(f)dx=
$$
$$
\int_{\mathbb{R}^n}\Lambda_e(f)\Omega_e(f^{p-1})dx+\int_{\mathbb{R}^n}\Lambda_o(f)\Omega_o(f^{p-1})dx=
$$
$$
=\int_{\mathbb{R}^n}\Lambda_e(f^{p-1})\Omega_e(f)dx+\int_{\mathbb{R}^n}\Lambda_o(f^{p-1})\Omega_o(f)dx.
$$
\end{theorem}
{\bf (III)} We turn our attention to the smoothness of the solutions of Lieb's equation (\ref{eq1}). We recall some notations and definitions 
from \cite{Vainikko}. Let $G\subseteq\mathbb{R}^n$ be an open and bounded set. For a $\lambda\in\mathbb{R}^n$, G. Vainniko introduces a weight
function:
$$
w_{\lambda}(x)=\left\{\begin{array}{lll} 1 & {\rm for} & \lambda<0 \\ (1+|\log\rho(x)|)^{-1} & {\rm for} & \lambda=0 \\
\rho(x)^{\lambda} & {\rm for} & \lambda>0 \end{array}\right.,\,\,\,x\in G,
$$
where $\rho(x)=\inf_{y\in\partial G} |x-y|$ is the distance from $x$ to the boundary $\partial G$ of $G$. Let $m\in\mathbb{Z}^+\cup\{0\}$,
$\nu\in\mathbb{R}$ satisfy $\nu<n$. We define the space $C^{m,\nu}(G)$ as the collection of all $m$ times continuously differentiable
functions $u:\,G\rightarrow\mathbb{R}$ (or $\mathbb{C}$) such that:
$$
||u||_{m,\nu}=\sum_{|\alpha|\le m}\sup_{x\in G}(w_{|\alpha|-(n-\nu)}(x)|D_{\alpha} u(x)|)<\infty.
$$
So $C^{m,\nu}(G)$ contains all the $m$ times continuously differentiable functions $u$ on $G$ whose derivatives near the boundary
$\partial G$ can be estimated as follows:
$$
|D_{\alpha}u(x)|\le\,{\rm Const.}\left\{\begin{array}{lll} 1 & {\rm for} & |\alpha|<n-\nu \\ 1+|\log\rho(x)| & {\rm for} & |\alpha|=n-\nu \\
\rho(x)^{n-\nu-|\alpha|} & {\rm for} & |\alpha|>n-\nu \end{array}\right.,\,\,\,x\in G,\,\,|\alpha|\le m.
$$
\begin{remark}
1) The function $||\cdot ||_{m,\nu}$ on the space $C^{m,\nu}(G)$ is a norm. \\
2) The space $(C^{m,\nu}(G),||\cdot ||_{m,\nu})$ is complete, i.e. it is a Banach space.
\end{remark}
Consider the following integral equation:
\begin{equation}
\label{eq11}
u(x)=\int_G K(x,y,u(y))dy+f(x),\,\,\,x\in G.
\end{equation}
we assume that the kernel $K(x,y,u)$ is $m$ times ($m\ge 1$) continuously differentiable with respect to $x,y,u$ for $x\in G$, $y\in G$,
$x\ne y$, $u\in\mathbb{R}$. We also assume that there is a real number $\nu<n$, such that, $\forall\,k\in\mathbb{Z}^+$ and $\alpha,\beta\in(\mathbb{Z}^+)^n$
for which $k+|\alpha|+|\beta|\le m$, the following inequalities hold:
\begin{equation}
\label{eq12}
\end{equation}
$$
D_{\alpha}^x D_{\beta}^{x+y}\frac{\partial^k}{\partial u^k}K(x,y,u)\le b_1(u)\left\{\begin{array}{lll} 1 & {\rm for} & \nu+|\alpha |<0 \\
1+|\log |x-y|| & {\rm for} & \nu+|\alpha |=0 \\ |x-y|^{-\nu-|\alpha |} & {\rm for} & \nu+|\alpha |>0 \end{array}\right.,
$$
\begin{equation}
\label{eq13}
\end{equation}
$$
|D_{\alpha}^x D_{\beta}^{x+y}\frac{\partial^k}{\partial u^k}K(x,y,u_1)-D_{\alpha}^x D_{\beta}^{x+y}\frac{\partial^k}{\partial u^k}K(x,y,u_2)|\le
$$
$$
\le b_2(u_1,u_2)|u_1-u_2|\left\{\begin{array}{lll} 1 & {\rm for} & \nu+|\alpha |<0 \\
1+|\log |x-y|| & {\rm for} & \nu+|\alpha |=0 \\ |x-y|^{-\nu-|\alpha |} & {\rm for} & \nu+|\alpha |>0 \end{array}\right..
$$
The functions $b_1:\,\mathbb{R}\rightarrow\mathbb{R}^+$ and $b_2:\,\mathbb{R}^2\rightarrow\mathbb{R}^+$ are assumed to be bounded on every
bounded region of $\mathbb{R}$ and $\mathbb{R}^2$ respectively. The notation of G. Vainikko reads as follows:
\begin{equation}
\label{eq14}
\end{equation}
$$
D_{\alpha}^x D_{\beta}^{x+y}=\left(\frac{\partial}{\partial x_1}\right)^{\alpha_1}\ldots\left(\frac{\partial}{\partial x_n}\right)^{\alpha_n}
\left(\frac{\partial}{\partial x_1}+\frac{\partial}{\partial y_1}\right)^{\beta_1}\ldots \left(\frac{\partial}{\partial x_n}+\frac{\partial}{\partial y_n}
\right)^{\beta_n}. 
$$
We can now conveniently quote the result we need from \cite{Vainikko}: \\
\\
{\bf Theorem 8.1.(\cite{Vainikko})} {\it Let $f\in C^{m,\nu}(G)$ (in (\ref{eq11})) and let the kernel $K(x,y,u)$ satisfy inequalities (\ref{eq12})
and (\ref{eq13}). If the integral equation (\ref{eq11}) has a solution $u\in L^{\infty}(G)$, then $u\in C^{m,\nu}(G)$.} \\
\\
\begin{remark}
1) There is a companion result (Theorem 8.2) in \cite{Vainikko}, but we will not use it here. \\
2) We note that if we restrict Lieb's integral equation (\ref{eq1}) to a bounded open $G\subseteq\mathbb{R}^n$, then it satisfies the assumptions
of Theorem 8.1. in \cite{Vainikko}. In this case $f(x)\equiv 0\in C^{\infty}(G)$ and $K(x,y,u)=|x-y|^{-\lambda}$ is independent of $u$. Also we note
that:
$$
\left(\frac{\partial}{\partial x_i}+\frac{\partial}{\partial y_i}\right)|x-y|^{-\lambda}\equiv 0,
$$
and so in inequalities (\ref{eq12}) and (\ref{eq13}) the only interesting values of $(k,|\alpha|,|\beta|)$ are $k=0$, $|\beta|=0$, $|\alpha|\le m$ and 
within the inequalities we start with $\nu+|\alpha|=\lambda$. Thus we can make effective calculations of the smoothness degree of a bounded solution
of the restriction to $G$ of Leib's integral equation (\ref{eq1}). \\
3) We note that any solution $f(x)$ of the Lieb integral equation (\ref{eq1}), must decay to zero at infinity, in order for the integral
$$
\int_{\mathbb{R}^n}|x-y|^{-\lambda}f(y)dy,\,\,\,\,0<\lambda<n,
$$
to converge at infinity. Thus any solution is bounded outside a ball $B_n(R)=\{x\in\mathbb{R}^n\,|\,|x|<R\}$ for a large enough radius $R$. Thus 
such a solution can have singularities only within the ball $B_n(R)$, and the solution tends to infinity at each such a singularity, otherwise by Theorem 8.1 in \cite{Vainikko} it will be smooth (of some degree) at such a point.
\end{remark}
We thus have the following:

\begin{theorem}
Let $f(x)$ be a solution of Lieb's integral equation (\ref{eq1}). Then $\lim_{|x|\rightarrow\infty}f(x)=0$, and there is a ball $B_n(R)$ such that $f(x)$
is bounded for $x\not\in B_n(R)$, $f(x)$ can have a finite set of singularities inside $B_n(R)$, and $f(x)$ tends to infinity when $x\rightarrow s$ for
each singular point $s$ of $f(x)$.
\end{theorem}
$\qed$

\section{Lieb's integral equation has at least two non-equivalent symmetric decreasing solutions, Theorem 1.1}

{\bf A proof.} We use the following formula for the Fourier transform, \cite{SteinWeiss}:
$$
\widehat{|y|^{-\nu}}=\int_{\mathbb{R}^n}|y|^{-\nu}\exp(-2\pi ix\cdot y)dy=\left\{\pi^{\nu-(n/2)}\Gamma\left(\frac{n}{2}-\frac{\nu}{2}\right)\left/
\Gamma\left(\frac{\nu}{2}\right)\right.\right\}|x|^{\nu-n},
$$
where $0<\nu<n$. We compute the Fourier transform of our integral:
$$
\widehat{\left(\int_{\mathbb{R}^n}|t-y|^{-\lambda}f(y)dy\right)(x)}=\widehat{\left(\int_{\mathbb{R}^n}|t-y|^{-\lambda}\left(C(n,\lambda)|y|^{-(n-\lambda/2)}\right)dy\right)(x)}=
$$
$$
=C(n,\lambda)\left(\int_{\mathbb{R}^n}|y|^{-\lambda}\exp\left(ix\cdot y\right)dy\right)\left(\int_{\mathbb{R}^n}|y|^{-(n-\lambda/2)}
\exp\left(ix\cdot y\right)dy\right)=
$$
$$
=C(n,\lambda)\left\{\pi^{\lambda-n/2}\Gamma\left(\frac{n}{2}-\frac{\lambda}{2}\right)\left/\Gamma\left(\frac{\lambda}{2}\right)\right.\right\}
|x|^{\lambda-n}\times
$$
$$
\times\left\{\pi^{(n-\lambda/2)-n/2}\Gamma\left(\frac{n}{2}-\frac{(n-\lambda)}{2}\right)\left/\Gamma\left(\frac{n-\lambda}{2}\right)\right.\right\}
|x|^{n-(\lambda/2)-n}=
$$
$$
=C(n,\lambda)\left\{\pi^{n/2}\left(\Gamma\left(\frac{n}{2}-\frac{\lambda}{2}\right)\Gamma\left(\frac{\lambda}{4}\right)^2\right)\left/
\left(\Gamma\left(\frac{\lambda}{2}\right)\Gamma\left(\frac{n}{2}-\frac{\lambda}{4}\right)^2\right)\right.\right\}\widehat{|y|^{-\lambda/2}}(x)=
$$
$$
=C(n,\lambda)^{\lambda/(2n-\lambda)}\left(\widehat{|y|^{-(n-\lambda/2)(\lambda/(2n-\lambda)}}\right)(x)=
$$
$$
=\left(\widehat{\left(C(n,\lambda)|y|^{-(n-\lambda/2)}\right)^{\lambda/(2n-\lambda)}}\right)(x)=
$$
$$
=\left(\widehat{f\left(y\right)^{\lambda/(2n-\lambda)}}\right)(x)=\left(\widehat{f\left(y\right)^{p-1}}\right)(x).
$$
Hence:
$$
\int_{\mathbb{R}^n}|x-y|^{-\lambda}f(y)dy=f(x)^{p-1}.
$$
$\qed $ \\
\\

\section{Integral relations between two solutions of equation (\ref{eq1}), orthogonality and commutativity, Corollary 1.3, Theorem 1.5 and
Theorem 1.7}
In this section we are interested in the form of integral relations between two solutions of the Lieb integral equations. Thus we will
not bother with convergence issues and just formally expand the formulas in Corollary 1.3, in Theorem 1.5 and in Theorem 1.7. We note 
that the identity in Corollary 1.3 follows by Theorem 1.1 and by the case $\alpha=\beta=\overline{0}$ in equation (\ref{eq6}) of Theorem 1.5.
Also Theorem 1.7 follows by Theorem 1.5. Hence we need to prove only Theorem 1.5: \\
We start with Lieb's integral equation, (\ref{eq1}) and perform on it a partial differentiation with respect to $x_j$. We formally
differentiate under the integral sign. The result is:
$$
\int_{\mathbb{R}^n}\frac{\partial}{\partial x_j}\left(\left|x-y\right|^{-\lambda}\right)f(y)dy=\frac{\partial f(x)^{p-1}}{\partial x_j}.
$$
We note that:
$$
\frac{\partial}{\partial x_j}\left(\left|x-y\right|^{-\lambda}\right)=-\frac{\partial}{\partial y_j}\left(\left|x-y\right|^{-\lambda}\right).
$$
Hence:
$$
-\int_{\mathbb{R}^n}\frac{\partial}{\partial y_j}\left(\left|x-y\right|^{-\lambda}\right)f(y)dy=\frac{\partial f(x)^{p-1}}{\partial x_j}.
$$
By integration by parts we deduce the following:
$$
\int_{\mathbb{R}^n}\left|x-y\right|^{-\lambda}\frac{\partial f(y)}{\partial y_j}dy=\frac{\partial f(x)^{p-1}}{\partial x_j}.
$$
We iterate this argument and obtain:
$$
\int_{\mathbb{R}^n}\left|x-y\right|^{-\lambda}\frac{\partial^2 f(y)}{\partial y_k\partial y_j}dy=\frac{\partial^2 f(x)^{p-1}}
{\partial x_k\partial x_j}.
$$
Now an inductive argument implies:
\begin{equation}
\label{eq15}
\int_{\mathbb{R}^n}\left|x-y\right|^{-\lambda}D_{\alpha}^{(y)}\left(f(y)\right)dy=D_{\alpha}^{(x)}\left(f(x)^{p-1}\right).
\end{equation}
Another outer induction gives, finally:
$$
\int_{\mathbb{R}^n}\left|x-y\right|^{-\lambda}\Lambda\left(f(y)\right)dy=\Lambda\left(f(x)^{p-1}\right).
$$
Next, let $g(x)$ be one more solution of equation (\ref{eq1}), i.e.:
$$
\int_{\mathbb{R}^n}\left|x-y\right|^{-\lambda}g(y)dy=g(x)^{p-1}.
$$
By what we have already done, we have:
\begin{equation}
\label{eq16}
\int_{\mathbb{R}^n}\left|x-y\right|^{-\lambda}D_{\beta}^{(y)}\left(g(y)\right)dy=D_{\beta}^{(x)}\left(g(x)^{p-1}\right).
\end{equation}
We now use the double integration technique. We multiply equation (\ref{eq15}) by $D_{\beta}^{(x)}(g(x))$ and integrate
$\int_{\mathbb{R}^n}\ldots dx$:
$$
\int_{\mathbb{R}^n}D_{\beta}^{(x)}\left(g(x)\right)\int_{\mathbb{R}^n}\left|x-y\right|^{-\lambda}D_{\alpha}^{(y)}\left(f(y)\right)dydx=
\int_{\mathbb{R}^n}D_{\beta}^{(x)}\left(g(x)\right)D_{\alpha}^{(x)}\left(f(x)^{p-1}\right)dx.
$$
Reversing the order of integration on the left hand side gives:
$$
\int_{\mathbb{R}^n}D_{\alpha}^{(y)}\left(f(y)\right)\int_{\mathbb{R}^n}\left|x-y\right|^{-\lambda}D_{\beta}^{(x)}\left(g(x)\right)dxdy=
\int_{\mathbb{R}^n}D_{\alpha}^{(y)}\left(f(y)\right)D_{\beta}^{(y)}\left(g(y)^{p-1}\right)dy.
$$
Changing on the right hand side the name of the variable from $y$ to $x$, gives:
$$
\int_{\mathbb{R}^n}D_{\beta}^{(x)}\left(g(x)\right)D_{\alpha}^{(x)}\left(f(x)^{p-1}\right)dx=\int_{\mathbb{R}^n}D_{\alpha}^{(x)}
\left(f(x)\right)D_{\beta}^{(x)}\left(g(x)^{p-1}\right)dx.
$$
This proves equation (\ref{eq6}). This is the commutativity part. We now prove orthogonality. We start with equation (\ref{eq6}) in the special
case $f(x)=g(x)$. We get:
$$
\int_{\mathbb{R}^n}D_{\beta}^{(x)}\left(f(x)\right)D_{\alpha}^{(x)}\left(f(x)^{p-1}\right)dx=\int_{\mathbb{R}^n}D_{\alpha}^{(x)}
\left(f(x)\right)D_{\beta}^{(x)}\left(f(x)^{p-1}\right)dx.
$$
Integration by parts gives:
$$
\int_{\mathbb{R}^n}D_{\alpha}^{(x)}
\left(f(x)\right)D_{\beta}^{(x)}\left(f(x)^{p-1}\right)dx=(-1)^{|\alpha|}\int_{\mathbb{R}^n}f(x)\cdot D_{\alpha+\beta}^{(x)}
\left(f(x)\right)dx,
$$
while
$$
\int_{\mathbb{R}^n}D_{\alpha}^{(x)}
\left(f(x)^{p-1}\right)D_{\beta}^{(x)}\left(f(x)\right)dx=(-1)^{|\beta|}\int_{\mathbb{R}^n}f(x)\cdot D_{\alpha+\beta}^{(x)}
\left(f(x)\right)dx.
$$
We proved equation (\ref{eq7}). Equation (\ref{eq8}) is a consequence of equation (\ref{eq6}) and equation (\ref{eq7}). This proves
Corollary 1.3, Theorem 1.5 and Theorem 1.7. $\qed $

\noindent
{\it Ronen Peretz \\
Department of Mathematics \\ Ben Gurion University of the Negev \\
Beer-Sheva , 84105 \\ Israel \\ E-mail: ronenp@math.bgu.ac.il} \\ 
 
\end{document}